# Notions of Input to Output Stability

Eduardo Sontag[*]  
Dept. of Mathematics  
Rutgers University  
New Brunswick, NJ 08903  
sontag@control.rutgers.edu

Yuan Wang[†]  
Dept. of Mathematics  
Florida Atlantic University  
Boca Raton, FL 33431  
ywang@math.fau.edu

## Abstract

This paper deals with concepts of output stability. Inspired in part by regulator theory, several variants are considered, which differ from each other in the requirements imposed upon transient behavior. The main results provide a comparison among the various notions, all of which specialize to input to state stability (ISS) when the output equals the complete state.

## 1 Introduction

This paper addresses questions of *output stability* for general finite-dimensional control systems

$$\dot{x}(t) = f(x(t), u(t)), \quad y(t) = h(x(t)). \tag{1}$$

(Technical assumptions on $f$, $h$, and admissible inputs, are described later.) Roughly, a system (1) is "output stable" if, for any initial state, the output $y(t)$ converges to zero as $t \to \infty$. Inputs $u$ may influence this stability in different ways; for instance, one may ask that $y(t) \to 0$ only for those inputs for which $u(t) \to 0$, or just that $y$ remains bounded whenever $u$ is bounded. Such behavior is of central interest in control theory. As an illustration, we will review below how regulation problems can be cast in these terms, letting $y(t)$ represent a quantity such as a tracking error. Another motivation for studying output stability arises in classical differential equations: "partial" asymptotic stability (cf. [26]) is nothing but the particular case of our study in which there are no inputs $u$ and the coordinates of $y$ are a subset of the coordinates of $x$ (that is to say, $h$ is a projection on a subspace of the state space $\mathbb{R}^n$). The notion of output stability is also related to that of "stability with respect to two measures", cf. [9].

The main starting point for our work is the observation that there are many *different* ways of making *mathematically precise* what one means by "$y(t) \to 0$ for every initial state" (and, when there are inputs, "provided that $u(t) \to 0$"). These different definitions need not result in equivalent notions; one must decide how uniform is the rate of convergence of $y(t)$ to zero, and precisely how the magnitudes of inputs and initial states affect convergence.

Indeed, our previous work on input-to-*state* stability (ISS, for short) was motivated in much the same way. The concept of ISS was originally introduced in [16] to address the problem when $y = x$. Major theoretical results were developed in [19] and [22], and applications to control design can be found in, among others, [4], [6], [7], [8], [12], [15], [22], and [25], as well as in the recent work [13] as a foundation for the formulation of robust tracking.

Actually, the original paper [16] had already introduced a notion of input/output stability (IOS), but the theoretical effort until now was almost exclusively directed towards the ISS special case. (There are two ways to formulate the property of input/output stability and its variants. One is in purely input/output terms, where one uses past inputs in order to represent initial conditions. Another is in state space terms, where the effect of past inputs is summarized by an initial state. In [16] an i/o approach was used, but here, because of our interest in initial-state dependence, we adopt the latter point of view. The relations between both approaches are explained in [16] and in more detail in [10] and [5].)

---

[*]This research was supported in part by US Air Force Grant F49620-98-1-0242  
[†]This research was supported in part by NSF Grant DMS-9457826

It turns out that the IOS case is substantially more complicated than ISS, in the sense that there are subtle possible differences in definitions. *One of the main objectives of this paper is to elaborate on these differences and to compare the various definitions*; the companion paper [24] provides Lyapunov-theoretic characterizations of each of them.

A second objective is to prove a theorem on *output redefinition* which (a) extends one of the main steps in linear regulation theory to general nonlinear systems, and (b) provides one of the main technical tools needed for the construction of Lyapunov functions in [24].

The organization of this paper is as follows. Section 1.1 starts with the review of certain facts from regulation theory; this material is provided merely as an additional motivation for our study, and is not required in order to follow the paper. Because of the technical character of the paper, it seems appropriate to provide an intuitive overview; thus, the rest of that section describes the main results in very informal terms. After that, in Section 2, we define our notions carefully and state precisely the main results. The rest of the paper contains the proofs. A preliminary version of this paper appeared in [23].

## 1.1 Additional Motivation and Informal Discussion

Output regulation problems encompass the main typical control objectives, namely, the analysis of feedback systems with the following property: for each exogenous signal $d(\cdot)$ (which might represent a disturbance to be rejected, or a signal to be tracked), the output $y(\cdot)$ (respectively, a quantity being stabilized, or the difference between a certain variable in the system and its desired target value) must decay to zero as $t \to \infty$. Typically (see e.g. [18], Section 8.2, or [14], Chapter 15 for linear systems, and [3], Chapter 8, for nonlinear generalizations), the exogenous signal is unknown but is constrained to lie in a certain prescribed class (for example, the class of all constant signals). Moreover, this class can be characterized through an "exosystem" given by differential equations (for example, the constant signals are precisely the possible solutions of $\dot{d} = 0$, for different initial conditions).

In order to focus on the questions of interest for this paper, we assume that we already have a closed-loop system exhibiting the desired regulation properties, ignoring the question of how an appropriate feedback system has been designed. Moreover, let us, for this introduction, restrict ourselves to linear time-invariant systems (local aspects of the theory can be generalized to certain nonlinear situations employing tools from center manifold theory, see [3]). The object of study becomes:

$$\begin{aligned} \dot{z} &= Az + Pw \\ \dot{w} &= Sw \\ y &= Cz + Qw, \end{aligned}$$

seen as a system $\dot{x} = f(x)$, $y = h(x)$, where the extended state $x$ consists of $z$ and $w$; the $z$-subsystem incorporates both the state of the system being regulated (the plant) and the state of the controller, and the equation $\dot{w} = Sw$ describes the exosystem that generates the disturbance or tracking signals of interest. This is a system without inputs; later we explain how inputs may be introduced into the model as well.

As an illustration, take the stabilization of the position $y$ of a second order system $\ddot{y} - y = u + w$ under the action of all possible constant disturbances $w$. The conventional proportional-integral-derivative (PID) controller uses a feedback law $u(t) = c_1 q(t) + c_2 y(t) + c_3 v(t)$, for appropriate gains $c_1, c_2, c_3$, where $q = \int y$ and $v = \dot{y}$. Let us take $c_1 = -1$, $c_2 = c_3 = -2$. If we view the disturbances as produced by the "exosystem" $\dot{w} = 0$, the complete system becomes

$$\dot{q} = y, \quad \dot{y} = v, \quad \dot{v} = -q - y - 2v + w, \quad \dot{w} = 0$$

with output $y$. That is, the plant/controller state $z$ is $\operatorname{col}(q, y, v)$, and

$$S = Q = 0, \; A = \begin{pmatrix} 0 & 1 & 0 \\ 0 & 0 & 1 \\ -1 & -1 & -2 \end{pmatrix}, \; C = \begin{pmatrix} 0 & 1 & 0 \end{pmatrix}, \; P = \begin{pmatrix} 0 \\ 0 \\ 1 \end{pmatrix}.$$

In linear regulator theory, the routine way to verify that the regulation objective has been met is as follows. Suppose that the matrix $A$ is Hurwitz and that there is some matrix $\Pi$ such that the following two identities ("Francis' equations") are satisfied:

$$\Pi S = A\Pi + P$$
$$0 = C\Pi + Q.$$

(The existence of $\Pi$ is necessary as well as sufficient for regulation, provided that the problem is appropriately posed, cf. [2] and [3].) Consider the new variable $\widehat{y} := z - \Pi w$. The first identity for $\Pi$ allows decoupling $\widehat{y}$ from $w$, leading to $\dot{\widehat{y}} = A\widehat{y}$. Since $A$ is a Hurwitz matrix, one concludes that $\widehat{y}(t) \to 0$ for all initial conditions. As the second identity for $\Pi$ gives that $y(t) = C\widehat{y}(t)$, one has the desired conclusion that $y(t) \to 0$.

Let us now express this convergence in a much more informative form. For that purpose, we introduce the map $\widehat{h} : x = (z, w) \mapsto |z - \Pi w| = |\widehat{y}|$. We also denote, for ease of future reference, $\chi(r) := r/|C|$ and $\beta(r,t) = r\left|e^{tA}\right|/|C|$, using $|\cdot|$ to denote Euclidean norm of vectors and also the corresponding induced matrix norm. So, $y = C\widehat{y}$ gives

$$\chi(|y|) = \chi(|h(x)|) \leq \widehat{h}(x) = |\widehat{y}| \tag{2}$$

for all $x$, and we also have $|y(t)| \leq \beta(|\widehat{y}(0)|, t)$ and in particular

$$|y(t)| \leq \beta(|x(0)|, t), \quad \forall t \geq 0 \tag{3}$$

along all solutions. This estimate quantifies the rate of decrease of $y$ to zero, and its overshoot, in terms of the initial state of the system. For the auxilliary variable $\widehat{y}$, we have in addition the following "stability" property:

$$|\widehat{y}(t)| \leq \sigma(|\widehat{y}(0)|), \quad \forall t \geq 0 \tag{4}$$

where $\sigma(r) := r \sup_{t \geq 0} \left|e^{tA}\right|$.

The use of $\widehat{y}$ (or equivalently, finding a solution $\Pi$ for the above matrix identities) is a key step in the analysis of regulation problems. Note the fundamental contrast between the behaviors of $\widehat{y}$ and $y$: because of (4), a zero initial value $\widehat{y}(0)$ implies $\widehat{y} \equiv 0$, which in regulation problems corresponds to the fact that the initial state of the "internal model" of the exosignal matches exactly the one for the exosignal; on the other hand, for the output $y$, typically an error signal, it may very well happen that $y(0) = 0$ but $y(t)$ is not identically zero. The fact about $\widehat{y}$ which allows deriving (3) is that $\widehat{y}$ dominates the original output $y$, in the sense of (2). *One of the main results in this paper provides an extension to very general nonlinear systems of the technique of output redefinition.*

To illustrate again with the PID example: one finds that $\Pi = \mathrm{col}\,(1,0,0)$ is the unique solution of the required equations, and the change of variables consists of replacing $q$ by $q - w$, the difference between the internal model of the disturbance and the disturbance itself, and $\widehat{h}(q, y, v, w) = (|q - w|^2 + |y|^2 + |v|^2)^{1/2}$. For instance, with $x(0) = y(0) = v(0) = 0$ and $w(0) = 1$ we obtain the output $y(t) = \frac{1}{2}t^2 e^{-t}$. Notice that this output has $y(0) = 0$ but is not identically zero, which is consistent with an estimate (3). On the other hand, the dominating output $\widehat{y} = (q - w, y, v))$ cannot exhibit such overshoot.

The discussion of regulation problems was for systems $\dot{x} = f(x)$ which are subject to no external inputs. This was done in order to simplify the presentation and because classically one does not consider external inputs. In general, however, one should study the effect on the feedback system of perturbations which were not exactly represented by the exosystem model. A special case would be, for instance, that in which the exosignals are not exactly modeled as produced by an exosystem, but have the form $w + u$, where $w$ is produced by an exosystem. Then one may ask if the feedback design is robust, in the sense that "small" $u$ implies a "small" asymptotic (steady-state) error for $y$, or that $u(t) \to 0$ implies $y(t) \to 0$. Experience with the notion of ISS then suggests that one should replace (3) by an estimate as follows:

$$\text{(IOS)} \qquad |y(t)| \;\leq\; \beta(|x(0)|, t) + \gamma(\|u\|).$$

By this we mean that for some functions $\gamma$ of class $\mathcal{K}$ and $\beta$ of class $\mathcal{KL}$ which depend only on the system being studied, and for each initial state and control, such an estimate holds for the ensuing output. We suppose as a standing hypothesis that the system is forward-complete, that is to say, solutions exist (and are unique) for $t \geq 0$, for any initial condition and any locally essentially bounded input $u$. (Recall that a function $\gamma : [0, \infty) \to [0, \infty)$ is of class $\mathcal{K}$ if it is strictly increasing and continuous, and satisfies $\gamma(0) = 0$, and of class $\mathcal{K}_\infty$ if it is also unbounded, and that $\mathcal{KL}$ is the class of functions $[0, \infty)^2 \to [0, \infty)$ which are of class $\mathcal{K}$ on the first argument and decrease to zero on the second argument.) The reason that inputs $u$ are not usually incorporated into the regulation problem statement is probably due to the fact that for *linear* systems it makes no difference: it is easy to see, from the variation of parameters formula, that IOS holds for all inputs if and only if it holds for the special case $u \equiv 0$.

Property (4) generalizes when there are inputs to the following "output Lagrange stability" property:

$$\text{(OL)} \qquad |y(t)| \;\leq\; \sigma_1(|y(0)|) + \sigma_2(\|u\|).$$

As mentioned earlier when discussing the classical tools of regulation theory, one of our main results is this: if a system satisfies IOS, then we can always find *another output*, let us call it $\widehat{y}$, which dominates $y$, in the sense of (2), and for which the estimate OL holds in addition to IOS.

As $\dot{\widehat{y}} = A\widehat{y}$ and $A$ is a Hurwitz matrix, in the case of linear regulator theory the redefined output $\widehat{y}$ satisfies a stronger decay condition, which in the input case leads naturally to an estimate as follows:

$$\text{(SIIOS)} \qquad |y(t)| \;\leq\; \beta(|y(0)|, t) + \gamma(\|u\|)$$

(we write $y$ instead of $\widehat{y}$ because we wish to define these notions for arbitrary systems). We will study this property as well. For linear systems, the conjunction of IOS and OL is equivalent to SIIOS. (Sketch of proof: with zero input and any initial state $x$ such that $Cx = 0$, OL gives us that $Ce^{tA}x \equiv 0$, which means that the kernel of $C$ coincides with the unobservable subspace $O(A, C)$. Therefore, referring once again to the Kalman observability decomposition and the notations in [18], Equation (6.8), now $y$ can be identified to the first $r$ coordinates of the state, which represent a stable system.) Remarkably, this equivalence breaks down for general nonlinear systems, as we will show.

Finally, as in the corresponding ISS paper [19], there are close relationships between output stability with respect to inputs, and robustness of stability under output feedback. This suggests the study of yet another property, which is obtained by a "small gain" argument from IOS: there must exist some $\chi \in \mathcal{K}_\infty$ so that

(ROS) $\qquad |y(t)| \leq \beta(|x(0)|, t) \quad \text{if } |u(t)| \leq \chi(|y(t)|) \quad \forall t.$

For linear systems, this property is equivalent to IOS, because applied when $u \equiv 0$ it coincides with IOS. One of our main contributions will be the construction of a counterexample to show that this equivalence also fails to generalize to nonlinear systems.

In summary, we will show that precisely these implications hold:

$$\text{SIIOS} \;\Rightarrow\; \text{OL \& IOS} \;\Rightarrow\; \text{IOS} \;\Rightarrow\; \text{ROS}$$

and show that under output redefinition the two middle properties coincide.

## 1.2 Related Notions

We caution the reader not to confuse IOS with the notion named *input/output to state stability* (IOSS) in [21] (also called "detectability" in [17], and "strong unboundedness observability" in [5]). This other notion roughly means that "no matter what the initial conditions are, if future inputs *and* outputs are small, the state must be eventually small". It is not a notion of stability; for instance, the unstable system $\dot{x} = x$, $y = x$ is IOSS. Rather, it represents a property of zero-state detectability. There is a fairly obvious connection between the various concepts introduced, however: a system is ISS if and only if it is both IOSS and IOS. This fact generalizes the linear systems theory result "internal stability is equivalent to detectability plus external stability" and its proof follows by routine arguments ([16, 10, 5]).

## 2 Definitions, Statements of Results

We assume, for the systems (1) being considered, that the maps $f : \mathbb{R}^n \times \mathbb{R}^m \to \mathbb{R}^n$ and $h : \mathbb{R}^n \to \mathbb{R}^p$ are locally Lipschitz continuous. We also assume that $f(0, 0) = 0$ and $h(0) = 0$. We use the symbol $|\cdot|$ for Euclidean norms in $\mathbb{R}^n$, $\mathbb{R}^m$, and $\mathbb{R}^p$.

By an *input* we mean a measurable and locally essentially bounded function $u : \mathcal{I} \to \mathbb{R}^m$, where $\mathcal{I}$ is a subinterval of $\mathbb{R}$ which contains the origin. Whenever the domain $\mathcal{I}$ of an input $u$ is not specified, it will be understood that $\mathcal{I} = \mathbb{R}_{\geq 0}$. The $L_\infty^m$-norm (possibly infinite) of an input $u$ is denoted by $\|u\|$, i.e. $\|u\| = (\text{ess}) \sup\{|u(t)|, t \in \mathcal{I}\}$. Given any input $u$ and any $\xi \in \mathbb{R}^n$, the unique maximal solution of the initial value problem $\dot{x} = f(x, u)$, $x(0) = \xi$ (defined on some maximal open subinterval of $\mathcal{I}$) is denoted by $x(\cdot, \xi, u)$. When $\mathcal{I} = \mathbb{R}_{\geq 0}$, this maximal subinterval has the form $[0, t_{\max})$, where $t_{\max} = t_{\max}(\xi, u)$. The system is said to be *forward complete* if for every initial state $\xi$ and for every input $u$ defined on $\mathbb{R}_{\geq 0}$, $t_{\max} = +\infty$. The corresponding output is denoted by $y(\cdot, \xi, u)$, that is, $y(t, \xi, u) = h(x(t, \xi, u))$ on the domain of definition of the solution.

## 2.1 Main Stability Concepts

**Definition 2.1** A forward-complete system is:

- *input to output stable* (IOS) if there exist a $\mathcal{KL}$-function $\beta$ and a $\mathcal{K}$-function $\gamma$ such that
$$|y(t,\xi,u)| \leq \beta(|\xi|,t) + \gamma(\|u\|), \qquad \forall\, t \geq 0; \tag{5}$$

- *output-Lagrange input to output stable* (OLIOS) if it is IOS and there exist some $\mathcal{K}$-functions $\sigma_1, \sigma_2$ such that
$$|y(t,\xi,u)| \leq \max\{\sigma_1(|h(\xi)|), \sigma_2(\|u\|)\}, \qquad \forall\, t \geq 0; \tag{6}$$

- *state-independent* IOS (SIIOS) there exist some $\beta \in \mathcal{KL}$ and some $\gamma \in \mathcal{K}$ such that
$$|y(t,\xi,u)| \leq \beta(|h(\xi)|,t) + \gamma(\|u\|), \qquad \forall\, t \geq 0. \tag{7}$$

In each case, we interpret the estimates as holding for all inputs $u$ and initial states $\xi \in \mathbb{R}^n$. □

To each given system (1) and each smooth function $\lambda : \mathbb{R}_{\geq 0} \to \mathbb{R}$, we associate the following system with inputs $d(\cdot)$:
$$\dot{x} = g(x,d) := f(x, d\lambda(|y|)), \quad y = h(x), \tag{8}$$

where $d \in \mathcal{M}_\mathcal{B}$, where $\mathcal{B}$ denotes the closed unit ball $\{|\mu| \leq 1\}$ in $\mathbb{R}^m$, and where, in general, $\mathcal{M}_\Omega$ denotes a system with controls restricted to take values in some subset $\Omega \subseteq \mathbb{R}^m$. For a given function $\lambda$, we let $x_\lambda(\cdot, \xi, d)$ (and $y_\lambda(\cdot, \xi, d)$, respectively) denote the trajectory (and the output function, respectively) of (8) corresponding to each initial state $\xi$ and each input $d(\cdot)$. Note that the system need not be complete even if the original system (1) is, but on the domain of definition of the solution, one has that $x_\lambda(t, \xi, d) = x(t, \xi, u)$, where $u(t) = \lambda(x_\lambda(t, \xi, d))$. Also note that, for system (8), $\xi = 0$ is an equilibrium for every $d$.

**Definition 2.2** A forward-complete system (1) is *robustly output stable* (ROS) if there exists a smooth $\mathcal{K}_\infty$-function $\lambda$ such that the corresponding system (8) is forward complete, and there exists some $\beta \in \mathcal{KL}$ such that
$$|y_\lambda(t,\xi,d)| \leq \beta(|\xi|,t), \tag{9}$$
for all $t \geq 0, \xi \in \mathbb{R}^n, d \in \mathcal{M}_\mathcal{B}$. □

The function $\lambda$ in (9) represents a robust output gain margin. It quantifies the magnitude of output feedback that can be tolerated by the system without destroying output stability.

**Remark 2.3** Suppose system (8) is forward complete. Then, the existence of a $\beta$ as in Property (9) is equivalent to the following:

1. there is a $\mathcal{K}_\infty$-function $\delta(\cdot)$ such that for any $\varepsilon > 0$, it holds that
$$|y_\lambda(t,\xi,d)| \leq \varepsilon, \qquad \forall\, t \geq 0, \ \forall\, d,$$
whenever $|\xi| \leq \delta(\varepsilon)$; and

2. for any $r \geq 0$ and any $\varepsilon > 0$, there exists some $T_{r,\varepsilon} > 0$ such that

$$|y_\lambda(t, \xi, d)| \leq \varepsilon$$

for all $t \geq T_{r,\varepsilon}$, all $d$, and all $|\xi| \leq r$.

This can be shown by following exactly the same steps as done in the proof of Proposition 2.5 of [11]. □

The following implications hold:

$$\text{SIIOS} \Rightarrow \text{OLIOS} \Rightarrow \text{IOS} \Rightarrow \text{ROS}$$

and all the reverse implications are false. The first two implications are clear from the definitions, and, not surprisingly, they cannot be reversed, as shown by counterexamples in Section 2.2. The last implication follows from this result, whose proof is provided in Section 3.2:

**Lemma 2.4** If a system is IOS, then it is ROS.

This lemma is a generalization of a "small gain" argument for ISS (for systems with full outputs $y = x$, IOS and ISS coincide), and its proof is analogous, although technically a bit more complicated, to the one given for the ISS case in [19]. Since in the special case of ISS the converse is true, it would be reasonable to expect that ROS and IOS be equivalent. Thus, the following result is surprising, cf. Section 4:

**Lemma 2.5** There is a ROS system which is not IOS.

As explained in the introduction, we will prove the following result. We say that a system (1) is OLIOS *under output redefinition* if there exist a locally Lipschitz map $h_0 : \mathbb{R}^n \to \mathbb{R}_{\geq 0}$ with $h(0) = 0$, and a $\chi \in \mathcal{K}_\infty$, such that

$$h_0(\xi) \geq \chi(|h(\xi)|)$$

for all $\xi$, so that the system

$$\dot{x} = f(x, u), \quad y = h_0(x) \tag{10}$$

is OLIOS.

**Theorem 2.1** *The following are equivalent for a system (1):*

1. *The system is IOS.*

2. *The system is OLIOS under output redefinition.*

## 2.2 Some Remarks

The estimate (5) can be restated in various equivalent ways. By causality, one may write $\gamma(\|u\|_{[0,t]})$ instead of $\gamma(\|u\|)$, where $\|u\|_{[0,t]}$ denotes the $L_\infty$ norm of $u$ restricted to the interval $[0, t]$. Similarly, the sum $\beta(|\xi|, t) + \gamma(\|u\|)$ could be replaced by $\max\{\beta(|\xi|, t), \gamma(\|u\|)\}$ (just use $2\beta$ and $2\gamma$).

When $h$ is the identity, IOS coincides with the by know well-known ISS property.

It is an easy exercise to show that for any class-$\mathcal{KL}$ function $\beta$ and any class-$\mathcal{K}$ function $\gamma$ there are a function $\widehat{\beta}$ and a class-$\mathcal{K}$ function $\widehat{\gamma}$ so that $\min\{\sigma(s), \beta(r,t)\} \leq \widehat{\beta}\left(s, \frac{t}{1+\widehat{\gamma}(r)}\right)$ for all $s, r, t$. It follows that a system is OLIOS if and only if there exist $\beta \in \mathcal{KL}$, $\rho \in \mathcal{K}$, and $\gamma \in \mathcal{K}$ such that

$$|y(t,\xi,u)| \leq \beta\left(|h(\xi)|, \frac{t}{1+\rho(|\xi|)}\right) + \gamma(\|u\|), \qquad \forall\, t \geq 0,$$

holds for all trajectories of the system. Note that this property nicely encapsulates both the IOS and output-Lagrange aspects of the OLIOS notion.

The paper [24] provides necessary and sufficient characterizations for each of the properties studied here, for systems that are uniformly bounded-input bounded-state (UBIBS), i.e., systems for which there exists some $\mathcal{K}$-function $\sigma$ such that the solution $x(t,\xi,u)$ is defined for all $t \geq 0$ and $|x(t,\xi,u)| \leq \max\{\sigma(|\xi|), \sigma(\|u\|)\}$ for all $t \geq 0$, every input $u$, and every initial state $\xi$. Among many other results, an equivalence is established between IOS and the existence of some smooth function $V : \mathbb{R}^n \to \mathbb{R}_{\geq 0}$ such that, for some $\alpha_1, \alpha_2 \in \mathcal{K}_\infty$, $\alpha_1(|h(\xi)|) \leq V(\xi) \leq \alpha_2(|\xi|)$ for all states $\xi$ and, for some $\chi \in \mathcal{K}$, $DV(\xi)f(\xi,\mu) < 0$ whenever $V(\xi) \geq \chi(|\mu|)$ (for all states $\xi$ and control values $\mu$). Such a $V$ can be alternatively interpreted as a new output $\widehat{y}$, which dominates the original output $y$, and which has the property that $\widehat{y}(t)$ converges to zero monotonically (at least while it is larger than a certain function of the current input). One obtains in this manner another output redefinition result. This appears to be the best possible result, since it follows from a counterexample given in [24] that, for systems with inputs, it is in general impossible to make a system SIIOS under output redefinition. (For systems without inputs, a redefinition result does hold.)

As promised, here are simple counterexamples to the first two implications. Consider the following system:

$$\dot{x}_1 = -x_1 + x_2 + u, \quad \dot{x}_2 = -x_1 - x_2 + u, \quad y = x_2. \tag{11}$$

System (11), being linear and stable, is ISS, and hence, it is IOS. However, this system is not OLIOS because $x_2(0) = 0$ and $u \equiv 0$ do not imply $x_2(t) \equiv 0$.

For an example of a OLIOS yet not SIIOS system, we consider:

$$\dot{x}_1 = 0, \quad \dot{x}_2 = -\frac{2x_2 + u}{1 + x_1^2}, \quad y = x_2, \tag{12}$$

System (12) is OLIOS. This can be shown by a Lyapunov approach. Let $V(\xi) = \xi_2^2/2$. Then, along any trajectory $x(t)$ of the system,

$$\dot{V}(x(t)) = -x_2(t)\frac{2x_2(t) + u(t)}{1 + x_1^2(t)} \leq -\frac{V(x(t))}{1 + x_1(0)^2}$$

whenever $\|u\| \leq \sqrt{V(x(t))}$. Using a standard comparison argument, (c.f. p441 of [16]), one can show that

$$V(x(t)) \leq \max\left\{V(x(0))e^{\frac{-t}{1+x_1(0)^2}}, \|u\|^2\right\}, \qquad \forall\, t \geq 0.$$

Consequently, the system is OLIOS. But this system is not SIIOS, as the decay rate of $x_2(t)$ depends on $x_1(0)$.

Finally, we remark that, for any linear system (not merely for those arising in regulation problems), the OLIOS output redefinition can be done in terms of linear functions. Indeed, suppose that the output $y = Cx$ of a linear system $\dot{x} = Ax$ satisfies $y(t) \to 0$ for all initial states, and let $\widehat{y}$ denote the states of the "observable part" of the system. Then $\widehat{y}$ is a stable variable and there is also a function $\chi$ so that (2) holds. (Sketch of proof, using the Kalman observability decomposition and the notations in [18], Equation (6.8): let $\widehat{y}$ denote the first $r$ coordinates of the state, so that $\dot{\widehat{y}} = A_1\widehat{y}$. By detectability of $(A_1, C_1)$, $y(t) \to 0$ implies that $\widehat{y}(t) \to 0$ for all solutions. Moreover, $y = C_1\widehat{y}$.)

## 3 Proofs

We will first prove Theorem 2.1, and then we will prove Lemma 2.4.

### 3.1 Proof of Theorem 2.1

An outline of the proof is as follows. It is natural, if we wish for future outputs to be small when the initial output is small and small inputs are applied, to start by defining a map $h_0(\xi)$ as the supremum, over all future inputs, of the difference $\|y\| - \gamma(\|u\|)$ (this is basically the value function for an $L^\infty$ differential game). After showing that the required estimates are satisfied, the next step is to show that $h_0$ is locally Lipschitz, and then to smooth-out $h_0$ away from the set where $h_0 = 0$. The final step is to "flatten" $h_0$ near the latter set.

Assume that system (10) with $y = h(x)$ is IOS. Thus, for all $\xi$ and $u$,
$$|y(t,\xi,u)| \leq \max\{\beta(|\xi|,t),\, \gamma(\|u\|)\}, \qquad \forall\, t \geq 0,$$
where $\beta \in \mathcal{KL}$ and, without loss of generality, $\gamma \in \mathcal{K}_\infty$. Let $h_0 : \mathbb{R}^n \to \mathbb{R}_{\geq 0}$ be defined by
$$h_0(\xi) = \sup_{t \geq 0, u} \{\max\{|y(t,\xi,u)| - \gamma(\|u\|),\, 0\}\}. \tag{13}$$

Observe that forward completeness is being used in this definition. Since $|y(0,\xi,u)| - \gamma(\|u\|) = |h(\xi)| \geq 0$ for $u \equiv 0$, the above is equivalent to
$$h_0(\xi) = \sup_{t \geq 0, u} \{|y(t,\xi,u)| - \gamma(\|u\|)\}.$$

It is clear that
$$|h(\xi)| \leq h_0(\xi) \leq \beta_0(|\xi|), \qquad \forall\, \xi \in \mathbb{R}^n,$$
where $\beta_0(s) = \beta(s,0)$. Since, for any $u_0$ with $\gamma(\|u_0\|) \geq \beta_0(s)$,
$$\max\{\beta(|\xi|,t),\, \gamma(\|u_0\|)\} - \gamma(\|u_0\|) \leq \max\{\beta_0(|\xi|),\, \gamma(\|u_0\|)\} - \gamma(\|u_0\|) = 0,$$
it follows that
$$h_0(\xi) = \sup_{t \geq 0, \|u\| \leq \gamma^{-1}(\beta_0(|\xi|))} \{|y(t,\xi,u)| - \gamma(\|u\|)\}. \tag{14}$$

Also note that for any $\tau \geq 0$ and any $v$,
$$\begin{aligned}
h_0(x(\tau,\xi,v)) &= \sup_{t \geq 0, u} \{|y(t, x(\tau,\xi,v), u)| - \gamma(\|u\|)\} \\
&= \sup_{t \geq 0, u} \{|y(t+\tau, \xi, v\sharp_\tau u)| - \gamma(\|v\sharp_\tau u\|_{[\tau,\infty)})\} \\
&\leq \sup_{t \geq 0, u} \{\beta(|\xi|, t+\tau) + \gamma(\|v\sharp_\tau u\|) - \gamma(\|v\sharp_\tau u\|_{[\tau,\infty)})\} \\
&\leq \beta(|\xi|, \tau) + \gamma(\|v\|_{[0,\tau)}),
\end{aligned} \tag{15}$$

where $v\sharp_\tau u$ is the concatenation of $v$ and $u$ defined by

$$v\sharp_\tau u(t) = \begin{cases} v(t), & \text{if } 0 \leq t < \tau, \\ u(t-\tau), & \text{if } t \geq \tau. \end{cases}$$

This shows that the system (10) with the output map $y = h_0(x)$ satisfies an IOS-type estimate (5), with the same functions $\beta$ and $\gamma$ as the original system.

Next, let us show that (10) with $y = h_0(x)$ also satisfies an output Lagrange estimate (6), with $\sigma_1(r) = 2r$ and $\sigma_2(r) = 2\gamma(r)$. Indeed, for any input $v$ and any $\tau \geq 0$, we have

$$\begin{aligned}
h_0(x(\tau,\xi,v)) &= \sup_{t\geq 0, u} \{|y(t, x(\tau,\xi,v), u)| - \gamma(\|u\|)\} \\
&= \sup_{t\geq 0, u} \{|y(t+\tau, \xi, v\sharp_\tau u)| - \gamma(\|v\sharp_\tau u\|_{[\tau,\infty)})\} \\
&\leq \sup_{s\geq 0, u} \{|y(s, \xi, v\sharp_\tau u)| - \gamma(\|v\sharp_\tau u\|) + \gamma(\|v\|)\} \\
&\leq \sup_{s\geq 0, \widetilde{u}} \{|y(s, \xi, \widetilde{u})| - \gamma(\|\widetilde{u}\|) + \gamma(\|v\|)\} \\
&= h_0(\xi) + \gamma(\|v\|) \leq \max\{2h_0(\xi),\ 2\gamma(\|v\|)\},
\end{aligned} \qquad (16)$$

as desired.

Define $C := \{\xi : h_0(\xi) = 0\}$. Then for any $\xi \notin C$, it holds that

$$h_0(\xi) = \sup_{0\leq t \leq t_\xi, \|u\| \leq \gamma^{-1}(\beta_0(|\xi|))} \{|y(t,\xi,u)| - \gamma(\|u\|)\},$$

where $t_\xi = T_{|\xi|}(h_0(\xi)/2)$, and $T_r(s)$ is associated with $\beta$ defined as in Lemma A.1.

**Lemma 3.1** The function $h_0$ is locally Lipschitz on the set where $h_0(\xi) \neq 0$ and continuous everywhere.

*Proof.* We first remark that

$$\varliminf_{\xi \to \xi_0} h_0(\xi) \geq h_0(\xi_0), \qquad \forall\, \xi_0 \in \mathbb{R}^n, \qquad (17)$$

that is, $h_0(\xi)$ is lower semi-continuous on $\mathbb{R}^n$. Indeed, pick $\xi_0$ and let $c := h(\xi_0)$. Take any $\varepsilon > 0$. Then there are some $u_0$ and $t_0$ so that $|y(t_0, \xi_0, u_0)| - \gamma(\|u_0\|) \geq c - \varepsilon/2$. By continuity of $y(t_0, \cdot, u_0)$, there is some neighborhood $\widetilde{U}_0$ of $\xi_0$ so that $|y(t_0, \xi, u_0)| - \gamma(\|u_0\|) \geq c - \varepsilon$ for all $\xi \in \widetilde{U}_0$. Thus $h_0(\xi) \geq c - \varepsilon$ for all $\xi \in \widetilde{U}_0$, and this establishes (17).

Fix any $\xi_0 \notin C$, and let $c_0 = h_0(\xi_0)/2$. Then there exists a bounded neighborhood $U_0$ of $\xi_0$ such that

$$h_0(\xi) \geq c_0, \qquad \forall\, \xi \in U_0.$$

Let $s_0$ be such that $|\xi| \leq s_0$ for all $\xi \in U_0$. Then

$$h_0(\xi) = \sup\{|y(t,\xi,u)| - \gamma(\|u\|) : t \in [0, t_1], \|u\| \leq b\}, \qquad \forall\, \xi \in U_0,$$

where $t_1 = T_{s_0}(c_0/2)$, and $b = \gamma^{-1}(\beta_0(s_0))$. By [11, Proposition 5.5], one knows that $x(t,\xi,u)$ is Lipschitz in $\xi$ uniformly on the set $\|u\| \leq b$, $\xi \in U_0$, and $t \in [0, t_1]$, and therefore, so is $y(t,\xi,u)$. Let $L_1$ be a constant such that

$$|y(t,\xi,u) - y(t,\eta,u)| \leq L_1 |\xi - \eta|, \qquad \forall\, \xi, \eta \in U_0,\ \forall\, 0 \leq t \leq t_1,\ \forall\, \|u\| \leq b.$$

For any $\varepsilon > 0$ and any $\xi \in U_0$, there exist some $t_{\xi,\varepsilon} \in [0, t_1]$ and some $u_{\xi,\varepsilon}$ such that

$$h_0(\xi) \leq |y(t_{\xi,\varepsilon}, \xi, u_{\xi,\varepsilon})| - \gamma(\|u_{\xi,\varepsilon}\|) + \varepsilon.$$

It then follows that, for any $\xi, \eta \in U_0$, for any $\varepsilon > 0$,

$$\begin{aligned} h_0(\xi) - h_0(\eta) &\leq |y(t_{\xi,\varepsilon}, \xi, u_{\xi,\varepsilon})| - \gamma(\|u_{\xi,\varepsilon}\|) + \varepsilon - (|y(t_{\xi,\varepsilon}, \eta, u_{\xi,\varepsilon})| - \gamma(\|u_{\xi,\varepsilon}\|)) \\ &\leq L_1 |\xi - \eta| + \varepsilon. \end{aligned}$$

Consequently,
$$h_0(\xi) - h_0(\eta) \leq L_1 |\xi - \eta|, \qquad \forall \xi, \eta \in U_0.$$

By symmetry,
$$h_0(\eta) - h_0(\xi) \leq L_1 |\xi - \eta|, \qquad \forall \xi, \eta \in U_0.$$

This proves that $h_0$ is locally Lipschitz on $\mathbb{R}^n \setminus C$.

We now show that $h_0$ is continuous on $C$. Fix $\xi_0 \in C$. One would like to show that

$$\lim_{\xi \to \xi_0} h_0(\xi) = 0. \tag{18}$$

Assume that this does not hold. Then there exists a sequence $\{\xi_k\}$ with $\xi_k \to \xi_0$ and some $\varepsilon_0 > 0$ such that $h_0(\xi_k) > \varepsilon_0$ for all $k$. Without loss of generality, one may assume that

$$|\xi_k| \leq s_1, \qquad \forall k,$$

for some $s_1 \geq 0$. It then follows that

$$h_0(\xi_k) = \sup \{|y(t, \xi_k, u)| - \gamma(\|u\|) : t \in [0, t_2], \|u\| \leq b_1\},$$

where $t_2 = T_{s_1}(\varepsilon_0/2)$, and $b_1 = \gamma^{-1}(\beta_0(s_1))$. Hence, for each $k$, there exists some $u_k$ with $\|u_k\| \leq b_1$ and some $\tau_k \in [0, t_2]$ such that

$$|y(\tau_k, \xi_k, u_k)| - \gamma(\|u_k\|) \geq h_0(\xi_k) - \varepsilon_0/2 \geq \varepsilon_0/2. \tag{19}$$

Again, by the locally Lipschitz continuity of the trajectories, one knows that there is some $L_2 > 0$ such that

$$|y(t, \xi_k, u) - y(t, \xi_0, u)| \leq L_2 |\xi_k - \xi_0|, \qquad \forall k \geq 0, \, \forall 0 \leq t \leq t_2, \, \forall \|u\| \leq b_1.$$

Hence,
$$|y(\tau_k, \xi_0, u_k)| - \gamma(\|u_k\|) \geq |y(\tau_k, \xi_k, u_k)| - \varepsilon_0/4 - \gamma(\|u_k\|) \geq \varepsilon_0/4$$

for $k$ large enough, contradicting the fact that $h_0(\xi) = 0$. This shows that (18) holds on $C$. ∎

Now we show how to modify $h_0$ to get an output function $\widetilde{h}$ that is locally Lipschitz everywhere so that system (10) with $\widetilde{h}$ is OLIOS.

We first pick a function $\overline{h}(\xi)$ that is smooth on $\mathbb{R}^n \setminus C$ with the property

$$\frac{h_0(\xi)}{2} \leq \overline{h}(\xi) \leq 2h_0(\xi), \qquad \forall \xi \in \mathbb{R}^n.$$

This can be done according to, e.g., Theorem B.1 in [11]. According to Lemma 4.3 in [11], there exists a $\mathcal{K}_\infty$-function $\rho$ such that $\rho \circ \overline{h}$ is smooth everywhere. Let $\widetilde{h} = \rho \circ \overline{h}$. Note then that

$$\rho(h_0(\xi)/2) \leq \widetilde{h}(\xi) \leq \rho(2h_0(\xi)). \tag{20}$$

Combining this with the fact that $h_0(\xi) \geq |h(\xi)|$, one sees that
$$\widetilde{h}(\xi) \geq \chi(|h(\xi)|), \qquad \forall \xi,$$
where $\chi(s) = \rho(s/2)$. Because of (15), one has
$$\widetilde{h}(x(t,\xi,u)) \leq \widetilde{\beta}(|\xi|,t) + \widetilde{\gamma}(\|u\|), \qquad \forall t \geq 0,$$
where $\widetilde{\beta}(s,t) = \chi(8\beta(s,t))$, and $\widetilde{\gamma}(s) = \chi(8\gamma(s))$, and because of (20) and (16), one has
$$\begin{aligned}\widetilde{h}(x(t,\xi,u)) &\leq \rho(2h_0(x(t,\xi,u))) \leq \max\{\rho(4h_0(\xi)),\, \rho(4\gamma(\|u\|))\} \\ &\leq \max\{\rho(8\rho^{-1}(\widetilde{h}(\xi))),\, \rho(4\gamma(\|u\|))\}, \qquad \forall t \geq 0,\end{aligned}$$
that is,
$$\widetilde{h}(x(t,\xi,u)) \leq \max\{\widetilde{\sigma}_1(\widetilde{h}(\xi)),\, \widetilde{\sigma}_2(\|u\|)\}, \qquad \forall t \geq 0,$$
for all $\xi$ and all $u$, where $\widetilde{\sigma}_1(s) = \rho(8\rho^{-1}(s))$ and $\widetilde{\sigma}_2(s) = \rho(4\gamma(s))$. We conclude that system (10) with the output function $y = \widetilde{h}(x)$ is OLIOS. ∎

### 3.2 Proof of Lemma 2.4

An outline of the proof is as follows. After providing a generalization of the "small-gain" result for ISS given in [19], we show that any OLIOS system is ROS. Then, we use an output redefinition in order to transform our IOS system into a OLIOS one, and finally apply the result to the transformed system.

#### 3.2.1 Output Lagrange Stability Under Small Output Feedback

**Lemma 3.2** Assume that the system (1) is forward complete and admits an output-Lagrange estimate as in (6). Then, there is a smooth $\mathcal{K}_\infty$ function $\lambda$ such that the resulting system (8) is forward complete and
$$\sigma_2\Big(|d(t)|\,\lambda(|y_\lambda(t,\xi,d)|)\Big) \leq \frac{1}{2}|h(\xi)| \tag{21}$$
holds a.e. on $[0,\infty)$.

*Proof.* Let $\sigma_1, \sigma_2$ be $\mathcal{K}$-functions such that (6) holds. Without loss of generality, we assume that both are in $\mathcal{K}_\infty$ and that $\sigma_1(s) \geq s$ for all $s \geq 0$. Hence, $\sigma_1^{-1}(s) \leq s$ for all $s \geq 0$. Let $\lambda$ be any smooth $\mathcal{K}_\infty$-function such that
$$\sigma_2(\lambda(s)) < \frac{1}{4}\sigma_1^{-1}(s), \qquad \forall\, s > 0.$$

Below we show that with such a choice of $\lambda$, the resulting system (8) satisfies the desired properties.

Recall that we use $x_\lambda(t,\xi,d)$ and $y_\lambda(t,\xi,d)$ denote the trajectory and the output for system (8) corresponding to the initial state $\xi$ and the input function $d$. To show that system (8) is complete, we first prove that (21) holds a.e. on the maximal interval of definition $[0, t_{\max})$ of the solution.

Pick any $\xi$ and any $d$, and use simply $x_\lambda(t)$ and $y_\lambda(t)$ to denote $x_\lambda(t,\xi,d)$ and $y_\lambda(t,\xi,d)$ respectively. To prove (21) on $[0, t_{\max})$, it is enough to show that
$$\sigma_2(\lambda(|y_\lambda(t)|)) \leq \frac{1}{2}|h(\xi)| \tag{22}$$

for such $t$.

*Case 1:* $h(\xi) \neq 0$. Since $\sigma_2(\lambda(|y_\lambda(0)|)) = \sigma_2(\lambda(|h(\xi)|)) < \frac{1}{4}\sigma_1^{-1}(|h(\xi)|) \leq \frac{1}{4}|h(\xi)|$, it follows that $\sigma_2(\lambda(|y_\lambda(t)|)) \leq \frac{1}{4}|h(\xi)|$ for $t$ small enough. Let

$$t_1 = \inf\left\{t \in (0, t_{\max}) : \sigma_2(\lambda(|y(t)|)) > \frac{1}{2}|h(\xi)|\right\}$$

with $t_1 = t_{\max}$ if the set is empty. Suppose by way of contradiction that $t_1 < t_{\max}$. Then (22) holds on $[0, t_1)$, and hence, (21) holds a.e. on $[0, t_1)$. Note that on $[0, t_{\max})$, $y_\lambda(t) = y(t, \xi, u)$ with $u(t) = d(t)\lambda(|y_\lambda(t)|)$. With (6), one sees that $|y_\lambda(t)| \leq \sigma_1(|h(\xi)|)$ for all $0 \leq t \leq t_1$, and in particular, $|y(t_1)| \leq \sigma_1(|h(\xi)|)$. Consequently,

$$\sigma_2(\lambda(|y(t_1)|)) \leq \frac{1}{4}\sigma_1^{-1}(|y(t_1)|) \leq \frac{1}{4}|h(\xi)|,$$

contradicting the definition of $t_1$. Thus, (22) holds for all $t \in [0, t_{\max})$.

*Case 2:* $h(\xi) = 0$. In this case, it is enough to show that $y_\lambda(t) = 0$ for all $t \in [0, t_{\max})$. Suppose this is not true. Then there exists some $\varepsilon > 0$ and some $t_2 \in (0, t_{\max})$ such that $|y_\lambda(t_2)| \geq \varepsilon$. Let $0 < \varepsilon_0 < \varepsilon$ be such that $\lambda^{-1}(\sigma_2^{-1}(\varepsilon_0)) < \varepsilon/2$. Then there is some $\tau \in (0, t_2)$ such that $|y_\lambda(\tau)| = \varepsilon_0$. Applying (22) proved for case 1 to the new initial state $\xi_1 := x_\lambda(\tau)$, one sees that

$$|y_\lambda(t)| \leq \lambda^{-1}\left(\sigma_2^{-1}\left(\frac{1}{2}|y_\lambda(\tau)|\right)\right) \leq \lambda^{-1}\left(\sigma_2^{-1}(\varepsilon_0)\right) < \varepsilon/2$$

for all $t \in [\tau, t_{\max})$, and in particular, $|y_\lambda(t_2)| < \varepsilon/2$, a contradiction. This shows that $y_\lambda(t) = 0$ for all $t \in [0, t_{\max})$.

We have showed that in both cases, (22) holds for all $t \in [0, t_{\max})$, which implies that, for any $\xi$ and any $d$, the function $u(t) := d(t)\lambda(|y_\lambda(t, \xi, d)|)$ remains essentially bounded on $[0, t_{\max})$. Suppose $t_{\max} < \infty$, then, by the forward completeness property of system (1), the trajectory $x_\lambda(t, \xi, d)$ (which is in fact $x(t, \xi, u)$ with $u(t) = d(t)\lambda(|y_\lambda(t, \xi, d)|)$) is bounded on $[0, t_{\max})$. This contradicts the maximality of $t_{\max}$. Therefore, $t_{\max} = \infty$ for every $\xi$ and every $d$. Consequently, (21) holds for all $t \in [0, \infty)$. ∎

### 3.2.2 Output Asymptotic Stability Under Small Output Feedback

**Lemma 3.3** *If a system (1) is* OLIOS, *then it is* ROS.

*Proof.* Estimates as in (5) and (6) hold. Without loss of generality, we assume that the functions appearing in the latter satisfy $\sigma_1(r) \geq r$ for all $r$ and $\sigma_2 \in \mathcal{K}_\infty$. Furthermore, redefining if necessary $\beta$ and $\gamma$, we replace (5) by a maximum, and use the same $\sigma_2$ function, so that the following estimate holds:

$$|y(t, \xi, u)| \leq \max\{\beta(|\xi|, t), \sigma_2(\|u\|)\}, \qquad \forall t \geq 0. \tag{23}$$

Let the function $\lambda$ be as in Lemma 3.2 so that system (8) is forward complete and (21) holds for almost all $t \geq 0$. By [1, Lemma 2.3], there exist some $\mathcal{K}$-functions $\varrho_1, \varrho_2$ and some $c \geq 0$ such that

$$|x_\lambda(t, \xi, d)| \leq \varrho_1(t) + \varrho_2(|\xi|) + c \tag{24}$$

for all $\xi$, all $d$, and all $t \geq 0$.

Below we will show that system (8) satisfies the two properties listed in Remark 2.3. First, by (21) and (6), one has:

$$|y_\lambda(t,\xi,d)| \leq \max\left\{\sigma_1(|h(\xi)|), \frac{1}{2}|h(\xi)|\right\} = \sigma_1(|h(\xi)|) \leq \widetilde{\sigma}(|\xi|), \qquad \forall\, t \geq 0,$$

where $\widetilde{\sigma}$ is any $\mathcal{K}$-function such that $\sigma_1(|h(\xi)|) \leq \widetilde{\sigma}(|\xi|)$ for all $\xi$. Property 1 follows readily. To prove Property 2, we first show the following:

*Claim.* For each $r > 0$, $s > 0$, there is some $T_{r,s} > 0$ such that

$$t \geq T_{r,s},\ |\xi| \leq r,\ |h(\xi)| \leq s \Rightarrow |y_\lambda(t,\xi,d)| \leq s/2. \tag{25}$$

To prove the claim, note that by (23) and (21), one has, for all $\xi$ as in (25),

$$|y_\lambda(t,\xi,d)| \leq \max\left\{\beta(|\xi|,t), \frac{|h(\xi)|}{2}\right\} \leq \max\left\{\beta(r,t), \frac{s}{2}\right\}, \qquad \forall\, t \geq 0.$$

Since $\beta \in \mathcal{KL}$, there is some $T_{r,s} > 0$ such that $\beta(r,t) \leq s/2$ for all $t \geq T_{r,s}$, and consequently,

$$|y(t,\xi,u)| \leq \frac{s}{2}, \qquad \forall\, t \geq T_{r,s}.$$

This $T_{r,s}$ satisfies the requirements of the claim.

Let $\kappa$ be a $\mathcal{K}$-function such that $|h(\xi)| \leq \kappa(|\xi|)$ for all $\xi$. Let $\varepsilon > 0$ be given. Pick any $\xi \neq 0$. Let $r = |\xi|$. Then $|h(\xi)| \leq \kappa(r)$. Let $l > 0$ be such that $2^{-l}\kappa(r) < \varepsilon$. Let $s_1 = \kappa(r)$ and $s_i = s_{i-1}/2$ for $i \geq 2$. By (25), there is some $T_{r,s_1} > 0$ such that

$$|y_\lambda(t,\xi,d)| \leq s_1/2, \qquad \forall\, t \geq T_{r,s_1},\ \forall\, d.$$

By (24), one has

$$|x_\lambda(T_{r,s_1},\xi,d)| \leq \varrho_1(T_{r,s_1}) + \varrho_2(r) + c := r_2, \qquad \forall\, d.$$

Applying (25) to $r_2$ and $s_2$, one sees that there is some $T_{r_2,s_2}$ such that the following holds:

$$|y_\lambda(t + T_{r,s_1},\xi,d)| \leq s_2/2, \qquad \forall\, t \geq T_{r_2,s_2},\ \forall\, d.$$

Inductively, letting $\widetilde{T}_k = \sum_{i=1}^{k} T_{r_i,s_i}$ (where $r_1 = r$) and applying (25) to $s_{k+1}$ and

$$r_{k+1} := \varrho_1\left(\widetilde{T}_k\right) + \varrho_2(r) + c,$$

one sees that there is some $T_{r_{k+1},s_{k+1}}$ such that

$$\left|y_\lambda(t + \widetilde{T}_k,\xi,d)\right| \leq \frac{s_k}{2}, \qquad \forall\, t \geq T_{r_k,s_k},\ \forall\, d.$$

Finally, we let $T = \widetilde{T}_1 + \widetilde{T}_2 + \ldots + \widetilde{T}_l$. Then, for any $t \geq T$ and any $d$,

$$|y_\lambda(t,\xi,d)| \leq \frac{s}{2^l} < \varepsilon.$$

Observe that in the above argument, $T$ only depends on $|\xi|$ and $\varepsilon$. Thus, the system satisfies Property 2 of Remark 2.3. Consequently, system (8) admits an estimate (9). ∎

### 3.2.3 Proof of Lemma 2.4

Lemma 2.4 follows easily from Theorem 2.1 and Lemma 3.3. Suppose system (1) is IOS. By Theorem 2.1, there is some locally Lipschitz function $h_0$ and some $\mathcal{K}_\infty$-function $\chi$ as in Theorem 2.1 such that system (10) with $y = h_0(x)$ as output is OLIOS. By Lemma 3.3, there is some smooth $\lambda_0 \in \mathcal{K}_\infty$ such that the system

$$\dot{x} = f(x, d\lambda_0(|h_0(x)|)), \ y = h_0(x), \qquad d \in \mathcal{M}_\mathcal{B},$$

is forward complete, and there is some $\beta_0 \in \mathcal{KL}$ such that

$$|h_0(x_{\lambda_0}(t, \xi, d))| \leq \beta_0(|\xi|, t), \qquad \forall\, t \geq 0,\ \forall\, d. \tag{26}$$

Recall that $\chi \in \mathcal{K}_\infty$ such that $\chi(|h(\xi)|) \leq h_0(\xi)$ for all $\xi$. Let $\lambda = \lambda_0 \circ \chi \in \mathcal{K}_\infty$. Then, $\lambda(|h(\xi)|) \leq \lambda_0(|h_0(\xi)|)$ for all $\xi$.

Consider the system

$$\dot{x} = f(x, d\lambda(|h(x)|)), \ y = h(x), \qquad d \in \mathcal{M}_\mathcal{B}. \tag{27}$$

Pick any $\xi$ and any $d$. Suppose the corresponding solution $x_\lambda(t, \xi, d)$ is defined on $[0, T)$ for some $T \leq \infty$. For $t \in [0, T)$, let

$$d_0(t) = \frac{\lambda(|h(x_\lambda(t, \xi, d))|)}{\lambda_0(|h_0(x_\lambda(t, \xi, d))|)} d(t)$$

if $h(x_\lambda(t, \xi, d)) \neq 0$, and $d_0(t) = 0$ otherwise. Extend $d_0$ to $[0, \infty)$ by letting $d_0(t) = 0$ if $t \geq T$ in case $T < \infty$. Observe that $d_0 \in \mathcal{M}_\mathcal{B}$, and by the uniqueness property,

$$x_\lambda(t, \xi, d) = x_{\lambda_0}(t, \xi, d_0), \qquad \forall\, t \in [0, T). \tag{28}$$

Suppose $T < \infty$. Then $x_{\lambda_0}(t, \xi, d_0)$ remains bounded on $[0, T)$, and hence, $x_\lambda(t, \xi, d)$ remains bounded on $[0, T)$. This contradicts the maximality of $T$. Consequently, $T = \infty$. This shows that system (27) is forward complete, and thus, (28) holds on $[0, \infty)$. Observe that

$$|h(x_\lambda(t, \xi, d))| \leq \chi^{-1}\left(h_0(x_{\lambda_0}(t, \xi, d_0))\right),$$

from which it follows by (26) that

$$|y_\lambda(t, \xi, d)| \leq \beta(|\xi|, t), \qquad \forall\, t \geq 0,\ \forall\, d,$$

where $\beta(s, t) := \chi^{-1}(\beta_0(s, t))$. This shows that system (1) is ROS. ∎

## 4 Example

In this section we show, by means of a counterexample, that ROS and IOS are not equivalent. We will then modify the example to get a system which is in addition bounded-input bounded-output stable (UBIBS), thus showing that even under this very strong stability assumption, ROS does not imply IOS.

Consider the following two-dimensional system:

$$\begin{aligned} \dot{x} &= \rho(|u| - 1 - |y|)x - y\sigma(x, y), \\ \dot{y} &= \rho(|u| - 1 - |y|)y + x\sigma(x, y), \end{aligned} \tag{29}$$

with output $y$ (we write $x$ and $y$ instead of $x_1$ and $x_2$), where $\rho$ is defined by

$$\rho(s) = \begin{cases} -1 & \text{if } s < -1, \\ s & \text{if } |s| \leq 1, \\ 1 & \text{if } s > 1, \end{cases}$$

and $\sigma(x,y)$ is defined by

$$\sigma(x,y) = \sigma_0(|y|-4) + (1-\sigma_0(|y|-4))\frac{\frac{1}{\pi}+|y|}{\max\{1,|x|\}},$$

and $\sigma_0$ is any smooth function such that $\sigma_0(s) = 1$ if $s \geq 0$, $\sigma(s) = 0$ if $s \leq -2$ and $0 < \sigma_0(s) \leq 1$ for $s \in (-2, 0)$.

First observe that the system is forward complete. This is because for the function $V(x,y) := (x^2+y^2)/2$, it holds that

$$DV(x,y)f(x,y,u) = \rho(|u|-1-|y|)(x^2+y^2) \leq 2\rho(|u|)V(x,y) \leq 2|u|\,V(x,y),$$

which implies that, along any trajectory $(x(t), y(t))$, one has $V(x(t), y(t)) \leq V(x(0), y(0))e^{2\|u\|t}$. This shows that all trajectories corresponding to (locally) bounded inputs are defined for all $t \geq 0$. Also note that $DV(x,y)f(x,y,u) = \rho(|u|-1-|y|)(x^2+y^2)$ implies that

$$|y| \geq |u| \Rightarrow DV(x,y)f(x,y,u) \leq -(x^2+y^2). \tag{30}$$

This means that every solution of the system (8), with $\lambda$ the identity function, satisfies the estimate $x(t)^2 + y(t)^2 \leq (x(0)^2 + y(0)^2)e^{-t}$. In particular, the output $y(t)$ satisfies a decay estimate (9), and in fact much more holds, since the origin of the system is globally asymptotically stable. Next, we show that property (5) fails to hold for the system. We do this by finding a bounded control (namely, $u \equiv 5$) and an initial state such that $y(t)$ is unbounded.

Pick the input function $u_0 \equiv 5$. For this input, the system satisfies the equations

$$\begin{aligned} \dot{x} &= \rho(4-|y|)x - y\sigma(x,y), \\ \dot{y} &= \rho(4-|y|)y + x\sigma(x,y). \end{aligned} \tag{31}$$

Using polar coordinates, on $\mathbb{R}^2 \setminus \{0\}$, the equations can also be written as

$$\begin{aligned} \dot{r} &= \rho(4-|y|)r, \\ \dot{\theta} &= \sigma_0(|y|-4) + (1-\sigma_0(|y|-4))\frac{\frac{1}{\pi}+|y|}{\max\{1,|x|\}}. \end{aligned}$$

Note that $\dot{r} > 0$ when $|y| < 4$ and $\dot{r} < 0$ when $|y| > 4$. Also note that $\dot{\theta} > 0$, and $\dot{\theta} = 1$ when $|y| \geq 4$.

Pick a point $z_0 = (x_0, y_0)$ with $x_0 > 0$ large enough and $y_0 = 4$, and consider the trajectory $\varsigma(t) = (x(t), y(t))$ with the initial state $\varsigma(0) = z_0$. Let $D_0$ denote the region $\{y > 4\}$. Since $\dot{y} > 0$ at $z_0$ and $\dot{\theta} = 1$ on $D_0$, there exists $0 < t_1 \leq \pi$ such that $\varsigma(t) \in D_0$ for all $t \in (0, t_1)$, and $y(t_1) = 4$. Let $p = (x(t_1), y(t_1))$, and for any point $a \in \mathbb{R}^2$, let $r_a = |a|$, and let $x_a$ ($y_a$, respectively) denote the $x$-coordinate ($y$-coordinate, respectively) of $a$. Since $\dot{r} \geq -r$ in $D_0$, one sees that $r_p \geq r_{z_0}e^{-t_1} \geq r_{z_0}e^{-\pi}$. Suppose $x_0$ is large enough such that $|x_p| > 5$ (this is possible because $r_p^2 = x_p^2 + 16 \geq r_{z_0}^2 e^{-2\pi}$). Since $\dot{y} > 0$ on any point where $x > 0$ and $y = 4$, it is impossible to have $x_p > 0$. Hence, $x_p < -5$ (cf. figure 1). Let $D_1 = \{(x,y) : x < -5,\ 2 <$ $y < 4\}$, and let $D_2 = \{(x, y) : x < -5,\ 0 < y < 2\}$. Since $\dot{y} < 0$ on the line $y = 4$ where $x < 0$, and in $D_1 \cup D_2$, $\dot{x} < 0$ and

$$\dot{y} \leq y + x\sigma_0(|y| - 4) - x(1 - \sigma_0(|y| - 4))\frac{\frac{1}{\pi} + y}{x}$$
$$\leq y + \max\{x, -(\frac{1}{\pi} + y)\} \leq -\frac{1}{\pi},$$

it follows that there exist $t_1 < t_2 < t_3 < \infty$ such that $\varsigma(t) \in D_1$ for $t \in (t_1, t_2)$, $\varsigma(t) \in D_2$ for $t \in (t_2, t_3)$, and $y(t_3) = 0$. Furthermore, on $D_2$, $\dot{r} = r$ and $\dot{y} = y + x\frac{\frac{1}{\pi}+y}{-x} = -1/\pi$. This implies that $t_3 - t_2 = 2\pi$, and hence, $r_s = r_q e^{2\pi}$, where $q = \varsigma(t_2), s = \varsigma(t_3)$. Combining this with the fact that $\dot{r} \geq 0$ in $D_1$, one sees that $r_s \geq r_p e^{2\pi} \geq r_{z_0} e^{\pi}$.

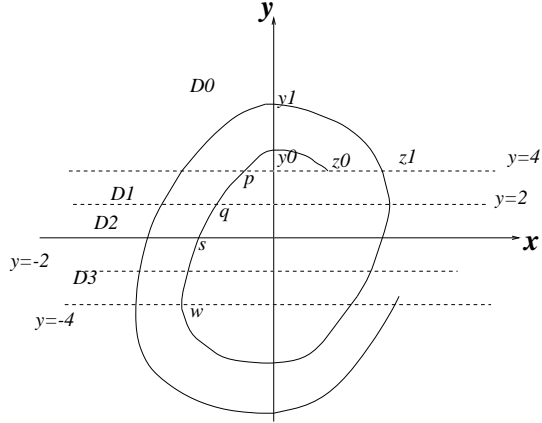

Figure 1: The trajectory of (31) starting from $z_0$

Let $D_3$ denote the region $\{(x, y) : x < -5, -4 < y < 0\}$. Then in $D_3$, $\dot{r} > 0$, $\dot{x} = \rho(4 - |y|)x - y\sigma(x, y) \leq |y|\max\{1, \frac{1}{\pi} + |y|\} \leq 20$, and,

$$\dot{y} \leq x\sigma_0(|y| - 4) + x(1 - \sigma_0(|y| - 4))\frac{\frac{1}{\pi} - y}{-x}$$
$$\leq \max\{x, -(\frac{1}{\pi} - y)\} \leq -\frac{1}{\pi},$$

it follows that if $|x_s|$ is large enough, then there exists $t_4 \in (t_3, t_3 + 4\pi]$ such that $\varsigma(t) \in D_3$ for all $t \in (t_3, t_4)$ and $y(t_4) = -4$. We assume that $x_0$ is large enough such that this happens (again, this is possible because $|x_s| \geq |x_p| = \sqrt{r_p^2 - 16}$ and $r_p^2 \geq (x_0^2 + 16)e^{-\pi}$). Let $w = \varsigma(t_4)$. Since in $D_3$ $\dot{r} \geq 0$, $r_w \geq r_s \geq r_{z_0} e^{\pi}$. Also observe that $|x_w| \geq x_0$ (because $r_w > r_{z_0}$, which implies $-x_w > x_0$, and $|y_w| = |y_0|$).

By symmetry, one can show that there exists $T_1 > 0$ such that $x(T_1) > 0$ and $y(T_1) = 4$ (i.e., at some moment $T_1$, $\varsigma(t)$ returns to the line $y = 4$), and moreover, with $z_1 = \varsigma(T_1), r_{z_1} \geq r_w e^{\pi} \geq r_{z_0} e^{2\pi}$.

**Remark 4.1** In the above we have seen that $\dot{\theta} = 1$ on $D_0$, and $\dot{y} \leq -1/\pi$ on $D_1 \cup D_2 \cup D_3$. Thus, $t_4 \leq t_1 + 8\pi \leq 9\pi$. By symmetry, one concludes that $T_1 \leq 18\pi$. □

Let $z_0, z_1, z_2, \ldots$ be the consecutive points where $\varsigma(t)$ intersects the line $y = 4$ where $x > 0$. The above argument shows that $r_{z_k} \geq r_{z_0} e^{2k\pi}$. Let $y_0, y_1, y_2, \ldots$ be the consecutive points where $\varsigma(t)$ intersects the $y$-axis in the upper half plane. Then (using once more that $\dot r \geq -r$ and $\dot\theta = 1$ on $D_0$) $y_k \geq r_{z_k} e^{-\pi/2} \geq r_{z_0} e^{2k\pi - \pi/2} \to \infty$. This shows that the output $y(t)$ corresponding to the initial state $z_0$ and the bounded input $u \equiv 5$ is unbounded, contradicting property (5).

Finally, we modify the above example in the following way to get a UBIBS system. Consider the two-dimensional system:

$$\begin{aligned}\dot x &= \rho(\varphi(x,y)\,|u| - 1 - |y|)x - y\sigma(x,y),\\ \dot y &= \rho(\varphi(x,y)\,|u| - 1 - |y|)y + x\sigma(x,y),\end{aligned} \quad (32)$$

where $\sigma, \rho$ are still defined the same as before, and the function $\varphi$ is a smooth function defined in the following way. For any point $\overline z = (\overline x, \overline y) \in C := \{(x,y) : x > 0, y = 4\}$ with $x$ large enough, let $T_{\overline z} = \inf\{t > 0 : \varsigma(t,\overline z) \in C\}$, where $\varsigma(t,\overline z)$ denotes the trajectory of (31) with $\varsigma(0,\overline z) = \overline z$. It was shown that $|\varsigma(T_z)| \geq |z|\,e^{2\pi}$. Now for each $k > 0$, let $A_k$ be the set $\{z \in \mathbb{R}^2 : r_k \leq |z| \leq r'_k\}$, where $0 < r_1 < r'_1 < r_2 < r'_2 < r_3 < r'_3 < \ldots \to \infty$ are such that, for any $k$, there exists some $z_k \in A_k \cap C$ so that for the trajectory $\varsigma(t, z_k)$ of (31), it holds that $\varsigma(t, z_k) \in A_k$ for all $t \in [0, T_{z_k}]$. Then $\varphi$ can be taken as any smooth function such that $\varphi(x,y) = 1$ on $A_k$, $\varphi(x,y) = 0$ for all $z = (x,y)$ such that $|z| = (r'_k + r_{k+1})/2$ for all $k \geq 1$, and $0 \leq \varphi(x,y) \leq 1$ everywhere else. In addition, we assume $\varphi(x,y) = 0$ on the set $A_0 := \{|z| \leq r_1/2\}$.

To show that the system is UBIBS, we still let $V(x,y) = (x^2 + y^2)/2$. Pick any input $u$ and any initial point $z_0 = (x_0, y_0)$. Let $(x(t), y(t))$ denote the corresponding trajectory. Let $\widetilde A_k$ be the set $\{|z| \leq r''_k := (r'_k + r_{k+1})/2\}$. Then $\widetilde A_k$ is forward invariant because on the set $\{|z_k| = r''_k\}$, $DV(x,y)f(x,y,u) = 2\rho(-1-|y|)V(x,y) < 0$. Hence, if $|z_0| \leq r''_k$, then $|(x(t), y(t))| \leq r''_k$ for all $t \geq 0$. Observe that $A_0$ is also forward invariant, and if $(x(0), y(0)) \in A_0$, $\frac{d}{dt} V(x(t), y(t)) \leq 2\rho(-1 - y(t))V(x(t), y(t)) \leq 0$, which implies that $V(x(t), y(t)) \leq V(x(0), y(0))$ for all $t \geq 0$. Now let $\sigma$ be any $\mathcal{K}_\infty$-function such that $\sigma(s) \geq s$ for $0 \leq s \leq r''_0$, and $\sigma(s) \geq r''_{k+1}$ for $r''_k \leq s \leq r''_{k+1}$ for all $k \geq 0$. Then $V(x(t), y(t)) \leq \sigma(|z(0)|)$ for all $t \geq 0$. This shows that the system is UBIBS. It can also be seen that $V$ still satisfies (30) for system (32), and hence, arguing as before, system (32) is ROS.

To show that the system fails to be IOS, we again pick the input function $u_0 \equiv 5$. With this input, the system satisfies the equations

$$\begin{aligned}\dot x &= \rho(5\varphi(x,y) - 1 - |y|)x - y\sigma(x,y),\\ \dot y &= \rho(5\varphi(x,y) - 1 - |y|)y + x\sigma(x,y).\end{aligned} \quad (33)$$

Observe that, for each $k \geq 1$, equations in (33) are the same as in (31) on $A_k$. For each $z \in \mathbb{R}^2$, let $\vartheta(t,z)$ denote the trajectory of (33) with the initial state $\vartheta(0,z) = z$. Then, for each $k \geq 1$, $\vartheta(t, z_k) = \varsigma(t, z_k)$ for all $t \in [0, T_{z_k}]$. Let $D_k$ be the complement of the region enclosed by the curve $\{\vartheta(t, z_k) : 0 \leq t \leq T_{z_k}\}$ and the line segment $L_k$ between $z_k$ and $\vartheta(T_{z_k})$ on the line $y = 4$. Then $D_k$ is forward invariant by uniqueness and the fact that $\dot y > 0$ on $L_k$. This implies that for all $t > T_{z_k}$, $|\vartheta(t, z_k)| \geq r_k$.

For each $k \geq 1$, let $(x_k(t), y_k(t)) = \vartheta(t, z_k)$. Below we will show that, for each $k$, there exists $t_1 < t_2 < t_3 < \ldots \to \infty$ such that $y_k(t_l) \geq r_k$ for all $l \geq 1$. It then will follow that $\varlimsup_{t\to\infty} y_k(t) \geq r_k$. For this purpose, we consider the angular movement $\theta_k(t)$ of $\vartheta(t, z_k)$. It can be seen that $\theta_k(t)$ satisfies the equation

$$\dot\theta = \sigma_0(|y|-4) + (1 - \sigma_0(|y|-4))\frac{\frac{1}{\pi} + |y|}{\max\{1, |x|\}}.$$

Since $\widetilde A_k$ is forward invariant, it holds that $|\vartheta(t, z_k)| \leq r''_k$ for all $t \geq 0$, and hence,

$$\frac{d}{dt}\theta_k(t) \geq \min\left\{1, \frac{\frac{1}{\pi} + |y_k(t)|}{\max\{1, |x_k(t)|\}}\right\} \geq \min\left\{1, \frac{\frac{1}{\pi}}{\max\{1, r''_k\}}\right\}.$$

From this, we know that $\theta_k(t) \to \infty$. This shows that there exist $0 < t_1 < t_2 < t_3 < \ldots \to \infty$ such that $\theta_k(t_l) = 2l\pi + \pi/2$. Hence, $|y_k(t_l)| = |\vartheta(t_l, z_k)| \geq r_k$ for all $l \geq 1$. This shows that it is impossible to have some $\gamma \in \mathcal{K}$ such that $\varlimsup_{t\to\infty} |y_k(t)| \leq \gamma(\|u_0\|) = \gamma(5)$ for all $k$. We conclude that the system is not IOS.

## A  A Lemma Regarding $\mathcal{KL}$ Functions

**Lemma A.1** For any $\mathcal{KL}$-function $\beta$, there exists a family of mappings $\{T_r\}_{r\geq 0}$ such that

- for each fixed $r > 0$, $T_r : \mathbb{R}_{>0} \xrightarrow{\text{onto}} \mathbb{R}_{>0}$ is continuous and strictly decreasing, and $T_0(s) \equiv 0$;

- for each fixed $s > 0$, $T_r(s)$ is strictly increasing as $r$ increases, and is such that $\beta(r, T_r(s)) < s$, and consequently, $\beta(r, t) < s$ for all $t \geq T_r(s)$.

*Proof.* For each $r \geq 0$ and each $s > 0$, let $\widehat{T}_r(s) := \inf\{t : \beta(r, t) \leq s\}$. Then $\widehat{T}_r(s) < \infty$, for any $r, s > 0$, $\beta(r, \widehat{T}_r(s)) \leq s$ for all $r \geq 0$, all $s > 0$, and it satisfies

$$\widehat{T}_r(s_1) \geq \widehat{T}_r(s_2), \text{ if } s_1 \leq s_2, \text{ and } \widehat{T}_{r_1}(s) \leq \widehat{T}_{r_2}(s), \text{ if } r_1 \leq r_2.$$

Note also that $\widehat{T}_0(s) = 0$ for all $s > 0$. Following exactly the same steps as in the proof of [11, Lemma 3.1], one can modify $\widehat{T}_r(s)$ to obtain $\widetilde{T}_r(s)$ so that for each fixed $r \geq 0$, $\widetilde{T}_r(\cdot)$ is decreasing and continuous; and for each fixed $s$, $\widetilde{T}_{(\cdot)}(s)$ is increasing. Finally, one lets $T_r(s) = \widehat{T}_r(s) + \frac{r}{1+s}$. Then $T_r(s)$ satisfies all conditions required in Lemma A.1. ∎